\documentclass[12pt]{article}
\usepackage[left=1.90cm,right=1.90cm,top=2.00cm,bottom=2.00cm]{geometry}
\usepackage{amsfonts}
\usepackage{titlesec}
\usepackage{cite,color,xcolor}
\usepackage{amsmath}
\usepackage{amssymb}
\usepackage{times}
\usepackage{bm}
\usepackage[colorlinks,citecolor=blue,urlcolor=blue]{hyperref}
\usepackage[utf8]{inputenc}
\usepackage{mathtools}
\usepackage{amsthm}
\usepackage{setspace}
\usepackage[numbers]{natbib}
\usepackage{cases}
\usepackage{fancyhdr}
\usepackage[pagewise]{lineno}

\newtheorem{theorem}{Theorem}[section]

\newtheorem{lemma}{Lemma}[section]

\newtheorem{definition}{Definition}[section]
\newtheorem{remark}{Remark}[section]

\usepackage{txfonts}

\newcommand{\bal}{\begin{align}}
\newcommand{\bbal}{\begin{align*}}
\newcommand{\bca}{\begin{cases}}
\newcommand{\eca}{\end{cases}}

\newcommand{\pa}{\partial}
\newcommand{\fr}{\frac}

\newcommand{\De}{\Delta}

\newcommand{\cd}{\cdot}

\newcommand{\dd}{\mathrm{d}}

\newcommand{\R}{\mathbb{R}}

\newcommand{\bi}{\Big}
\newcommand{\g}{\big}

\begin{document}
\title{Ill-posedness for the Camassa-Holm equation in $B_{p,1}^{1}\cap C^{0,1}$}

\author{Jinlu Li$^{1}$, Yanghai Yu$^{2,}$\footnote{E-mail: lijinlu@gnnu.edu.cn; yuyanghai214@sina.com(Corresponding author); guoyy35@fosu.edu.cn; mathzwp2010@163.com}, Yingying Guo$^{3}$ and Weipeng Zhu$^{3}$\\
\small $^1$ School of Mathematics and Computer Sciences, Gannan Normal University, Ganzhou 341000, China\\
\small $^2$ School of Mathematics and Statistics, Anhui Normal University, Wuhu 241002, China\\
\small $^3$ School of Mathematics and Big Data, Foshan University, Foshan, Guangdong 528000, China}

\date{\today}

\maketitle\noindent{\hrulefill}

{\bf Abstract:} In this paper, we study the Cauchy problem for the Camassa-Holm equation on the real line. By presenting a new construction of initial data, we show that the solution map in the smaller space $B_{p,1}^{1}\cap C^{0,1}$ with $p\in(2,\infty]$ is discontinuous at origin. More precisely, $u_0\in B_{p,1}^{1}\cap C^{0,1}$ can guarantee that the Camassa-Holm equation has a unique local solution in $W^{1,p}\cap C^{0,1}$, however, this solution is instable and can have an inflation in $B_{p,1}^{1}\cap C^{0,1}$ for certain initial data.


\vskip0mm\noindent{\hrulefill}

\section{Introduction}

In this paper we are concerned with the Cauchy problem for the classical Camassa-Holm (CH) equation
\begin{equation}\label{0}
\begin{cases}
u_t-u_{xxt}+3uu_x=2u_xu_{xx}+uu_{xxx}, \\
u(x,t=0)=u_0(x).
\end{cases}
\end{equation}
Here $(x,t)\in \R\times\R^+$, the scalar function $u = u(t, x)$ stands for the fluid velocity at time $t\geq0$ in the $x$ direction.
We can transform the CH equation equivalently into the following transport type equation
\begin{equation}\label{CH}
\begin{cases}
\partial_tu+u\pa_xu=-\pa_x\g(1-\pa^2_x\g)^{-1}\bi(u^2+\fr12(\pa_xu)^2\bi), \\
u(0,x)=u_0(x).
\end{cases}
\end{equation}
The CH equation was firstly proposed in the context of hereditary symmetries studied in \cite{Fokas} and then was derived explicitly as a water wave equation by Camassa-Holm \cite{Camassa}. Many aspects of the mathematical beauty of the CH equation have been exposed over the last two decades. Particularly, (CH) is completely integrable \cite{Camassa,Constantin-P} with a bi-Hamiltonian structure \cite{Constantin-E,Fokas} and infinitely many conservation laws \cite{Camassa,Fokas}. Also, it admits exact peaked
soliton solutions (peakons) of the form $$u(x,t)=ce^{-|x-ct|},\;c>0,$$ which are orbitally stable \cite{Constantin.Strauss}. Another remarkable feature of the CH equation is the wave breaking phenomena: the solution remains bounded while its slope becomes unbounded in finite time \cite{Constantin,Escher2,Escher3}. It is worth mentioning that the peaked solitons present the characteristic for the travelling water waves of greatest height and largest amplitude and arise as solutions to the free-boundary problem for incompressible Euler equations over a flat bed, see Refs. \cite{Constantin-I,Escher4,Escher5,Toland} for the details.

Due to these interesting and remarkable features, the CH equation has attracted much attention as a class of integrable shallow water wave equations in recent twenty years. Its systematic mathematical study was initiated in a series of papers by Constantin and Escher, see \cite{Escher1,Escher2,Escher3,Escher4,Escher5}. After the CH equation was derived physically in the context of water waves, there are
a large amount of literatures devoted to studying the well-posedness of the Cauchy problem \eqref{0} (see
Molinet's survey \cite{Molinet}). Particularly, the continuous dependence is rather important when PDEs are used to model
phenomena in the natural world since measurements are always associated with errors. Next we recall some progresses in this field.

{\bf Well-posedness.} Li and Olver \cite{li2000} proved that the Cauchy problem \eqref{0} is locally well-posed with the initial data $u_0\in H^s(\R)$ with $s > 3/2$ (see also \cite{GB}). Danchin \cite{d1} proved
the local existence and uniqueness of strong solutions to \eqref{0} with initial data in $B^s_{p,r}$ for $s > \max\{1 + 1/p , 3/2\}$ with $p\in[1,\infty]$ and $r\in[1,\infty)$. Meanwhile, he \cite{d1} only obtained the continuity of the solution map of \eqref{0} with respect to the initial data in the space $\mathcal{C}([0, T ];B^{s'}_{p,r})$ with any $s'<s$. Li-Yin \cite{Li-Yin1} proved the continuity of the solution map of \eqref{0} with respect to the initial data in the space $\mathcal{C}([0, T];B^{s}_{p,r})$ with $r<\infty$. In particular, they \cite{Li-Yin1} proved that the solution map of \eqref{0} is weak continuous with respect to initial data $u_0\in B^s_{p,\infty}$. For the endpoints, Danchin \cite{d3} obtained the local well-posedness in the space $B^{3/2}_{2,1}$. Recently, Ye-Yin-Guo \cite{Ye} proved the uniqueness and continuous dependence of the Camassa-Holm type equations in critical Besov spaces $B^{1+1/p}_{p,1}$ with $p\in[1,\infty)$.

{\bf Ill-posedness.} When considering further continuous dependence, we proved the non-uniform dependence on initial data for \eqref{0} under both the framework of Besov spaces $B^s_{p,r}$ for $s>\max\big\{1+1/p, 3/2\big\}$ with $p\in[1,\infty], r\in[1,\infty)$ and $B^{3/2}_{2,1}$ in \cite{Li1, Li2} (see \cite{H-M,H-K,H-K-M} for earlier results in $H^s$). Danchin \cite{d3} obtained the ill-posedness of \eqref{0} in $B^{3/2}_{2,\infty}$ (the data-to-solution map is not continuous by using peakon solution). Byers \cite{Byers} proved that the Camassa-Holm equation is ill-posed in $H^s$ for $s < 3/2$ in the sense of norm inflation, which means that $H^{3/2}$ is the critical Sobolev space for well-posedness. Moreover, Guo-Liu-Molinet-Yin \cite{Guo-Yin} showed the ill-posedness of \eqref{0} in $B_{p,r}^{1+1/p}(\mathbb{R}\;\text{or}\; \mathbb{T})$ with $(p,r)\in[1,\infty]\times(1,\infty]$ (especially in $H^{3/2}$) by proving the norm inflation. In our recent paper \cite{Li22}, we established the ill-posedness for \eqref{0} in $B^s_{p,\infty}(\mathbb{R})$ by proving the solution map to the CH equation starting from $u_0$ is discontinuous at $t = 0$ in the metric of $B^s_{p,\infty}(\mathbb{R})$. Very recently, Guo-Ye-Yin \cite{Guo} obtained the ill-posedness for the CH equation in $B^{1}_{\infty,1}(\R)$ by proving the norm inflation.

Now we would like to point out the idea that motivate the current work, which is inspired by Wang \cite{Wang} who proved ill-posedness for the 3D Navier-Stokes equation in critical Besov spaces $\dot{B}^{-1}_{\infty,q}$ in the sense of norm inflation. For arbitrarily small $\delta>0$, Wang \cite{Wang} showed that the solution map $\delta u_{0} \rightarrow u$ of the 3D Navier-Stokes equation in critical Besov spaces $\dot{B}_{\infty, q}^{-1}$ with any $q \in[1,2]$ is discontinuous at origin. Koch-Tataru \cite{Koch2001} proved that the 3D Navier-Stokes equation is globally well-posed for small data in $\mathrm{BMO}^{-1}$. Due to the embedding $\dot{B}_{\infty, q}^{-1} \subset \mathrm{BMO}^{-1}(q \leqslant 2)$, then for sufficiently small $\delta>0, u_{0} \in \dot{B}_{\infty, q}^{-1}(q \leqslant 2)$ can guarantee that the 3D Navier-Stokes equation has a unique global solution in $\mathrm{BMO}^{-1}$. However, Wang \cite{Wang} obtained that this solution is instable in $\dot{B}_{\infty, q}^{-1}$ and the solution can have an inflation in $\dot{B}_{\infty,q}^{-1}$ by constructing a special initial data. For the CH equation, we should note that, the solution map in $W^{1,p}\cap C^{0,1}$ is continuous at origin since the following estimate
$$\|u\|_{W^{1,p}\cap C^{0,1}}\leq \|u_0\|_{W^{1,p}\cap C^{0,1}}\exp\left(C\int_0^t\|u\|_{C^{0,1}}\dd \tau\right).$$
Taking notice of the inclusion $B^{1}_{p,1}\hookrightarrow W^{1,p}$, it seems natural to conjecture that the solution map of the CH equation in $B_{p,1}^{1}\cap C^{0,1}$ is also continuous at origin. However, we shall show the following negative result in this paper.
\begin{theorem}\label{th1}
Let $p\in(2,\infty]$. For any $n\in \mathbb{Z}^+$ large enough there exists $u_0$ with
$$\|u_0\|_{B^1_{p,1}\cap C^{0,1}}\leq \frac{1}{\log\log n},$$
such that if we denote by $u\in \mathcal{C}([0,1];H^{3})$, the solution of the Camassa-Holm equation with initial data $u_0$, then
$$\|u(t_0)\|_{B^1_{p,1}}\geq {\log\log n},$$
with $t_0\in (0,\frac{1}{\log n}]$.
\end{theorem}
\begin{remark}
Theorem \ref{th1} demonstrates the ill-posedness of the Camassa-Holm equation in $B^1_{p,1}\cap C^{0,1}$ in the sense of that Hadamard (the norm inflation implies discontinuity with respect to the initial data).
\end{remark}
\begin{remark}
Compared with the recent result by Guo-Ye-Yin \cite{Guo} who showed ill-posedness of the Camassa-Holm equation in $B^1_{\infty,1}$, Theorem \ref{th1} seems more general owing to $B^1_{\infty,1}\hookrightarrow C^{0,1}$. We should also mention that our result is not just a little bit of improvement. The key point is that, our construction of initial data is different essentially from that in \cite{Guo} since their initial data is no longer valid when proving the norm inflation in $B^1_{p,1}$ if $p\neq\infty$.
\end{remark}
This paper is structured as follows. In Section \ref{sec2}, we list some notations and known results and recall some Lemmas which will be used in the sequel. In Section \ref{sec3}, we prove Theorem \ref{th1}.
\section{Preliminaries}\label{sec2}
\subsection{Notations}\label{subsec21}
We will use the following notations throughout this paper. 
\begin{itemize}
  \item $C$ stands for some positive constant which may vary from line to line. 
  \item The symbol $A\approx_s B$ means that $c(s)B\leq A\leq C(s)B$. 
  \item We shall call a ball $B(x_0,r)=\{x\in \R: |x-x_0|\leq R\}$ with $R>0$ and an annulus $\mathcal{C}(0,r_1,r_2)=\{x\in \R: 0<r_1\leq|x|\leq r_2\}$ with $0 <r_1 <r_2$.
  \item Given a Banach space $X$, we denote its norm by $\|\cdot\|_{X}$.
  \item We denote by $W^{1,p}$ the standard Sobolev space on $\R$ of $L^p$ functions whose derivative also belongs to $L^p$.
\item We will also define the Lipschitz space $C^{0,1}$ using the norm $\|f\|_{C^{0,1}}=\|f\|_{L^\infty}+\|\pa_xf\|_{L^\infty}$.
\item For $I\subset\R$, we denote by $\mathcal{C}(I;X)$ the set of continuous functions on $I$ with values in $X$. Sometimes we will denote $L^p(0,T;X)$ by $L_T^pX$.
\item Let us recall that for all $f\in \mathcal{S}'$, the Fourier transform $\widehat{f}$, is defined by
$$
\widehat{f}(\xi)=\int_{\R}e^{-ix\xi}f(x)\dd x \quad\text{for any}\; \xi\in\R.
$$
\end{itemize}

\subsection{Littlewood-Paley analysis}\label{subsec22}
Next, we will recall some facts about the Littlewood-Paley decomposition, the nonhomogeneous Besov spaces and their some useful properties (see \cite{B} for more details).

Let $\mathcal{B}:=\{\xi\in\mathbb{R}:|\xi|\leq 4/3\}$ and $\mathcal{C}:=\{\xi\in\mathbb{R}:3/4\leq|\xi|\leq 8/3\}.$
Choose a radial, non-negative, smooth function $\chi:\R\mapsto [0,1]$ such that it is supported in $\mathcal{B}$ and $\chi\equiv1$ for $|\xi|\leq3/4$. Setting $\varphi(\xi):=\chi(\xi/2)-\chi(\xi)$, then we deduce that $\varphi$ is supported in $\mathcal{C}$. Moreover,
\begin{eqnarray*}
\chi(\xi)+\sum_{j\geq0}\varphi(2^{-j}\xi)=1 \quad \mbox{ for any } \xi\in \R.
\end{eqnarray*}
We should emphasize that the fact $\varphi(\xi)\equiv 1$ for $4/3\leq |\xi|\leq 3/2$ will be used in the sequel.

For every $u\in \mathcal{S'}(\mathbb{R})$, the inhomogeneous dyadic blocks ${\Delta}_j$ are defined as follows
\begin{numcases}{\Delta_ju=}
0, & if $j\leq-2$;\nonumber\\
\chi(D)u=\mathcal{F}^{-1}(\chi \mathcal{F}u), & if $j=-1$;\nonumber\\
\varphi(2^{-j}D)u=\mathcal{F}^{-1}\g(\varphi(2^{-j}\cdot)\mathcal{F}u\g), & if $j\geq0$.\nonumber
\end{numcases}
In the inhomogeneous case, the following Littlewood-Paley decomposition makes sense
$$
u=\sum_{j\geq-1}{\Delta}_ju\quad \text{for any}\;u\in \mathcal{S'}(\mathbb{R}).
$$
\begin{definition}[see \cite{B}]
Let $s\in\mathbb{R}$ and $(p,r)\in[1, \infty]^2$. The nonhomogeneous Besov space $B^{s}_{p,r}(\R)$ is defined by
\begin{align*}
B^{s}_{p,r}(\R):=\Big\{f\in \mathcal{S}'(\R):\;\|f\|_{B^{s}_{p,r}(\mathbb{R})}<\infty\Big\},
\end{align*}
where
\begin{numcases}{\|f\|_{B^{s}_{p,r}(\mathbb{R})}=}
\left(\sum_{j\geq-1}2^{sjr}\|\Delta_jf\|^r_{L^p(\mathbb{R})}\right)^{\fr1r}, &if $1\leq r<\infty$,\nonumber\\
\sup_{j\geq-1}2^{sj}\|\Delta_jf\|_{L^p(\mathbb{R})}, &if $r=\infty$.\nonumber
\end{numcases}
\end{definition}
The following Bernstein's inequalities will be used in the sequel.
\begin{lemma}[see Lemma 2.1 in \cite{B}] \label{lem2.1} Let $\mathcal{B}$ be a Ball and $\mathcal{C}$ be an annulus. There exist constants $C>0$ such that for all $k\in \mathbb{N}\cup \{0\}$, any positive real number $\lambda$ and any function $f\in L^p$ with $1\leq p \leq q \leq \infty$, we have
\begin{align*}
&{\rm{supp}}\hat{f}\subset \lambda \mathcal{B}\;\Rightarrow\; \|\pa_x^kf\|_{L^q}\leq C^{k+1}\lambda^{k+(\frac{1}{p}-\frac{1}{q})}\|f\|_{L^p},  \\
&{\rm{supp}}\hat{f}\subset \lambda \mathcal{C}\;\Rightarrow\; C^{-k-1}\lambda^k\|f\|_{L^p} \leq \|\pa_x^kf\|_{L^p} \leq C^{k+1}\lambda^k\|f\|_{L^p}.
\end{align*}
\end{lemma}
\begin{lemma}[see Lemma 2.100 in \cite{B}] \label{lem2.2}
Let $1 \leq r \leq \infty$, $1 \leq p \leq p_{1} \leq \infty$ and $\frac{1}{p_{2}}=\frac{1}{p}-\frac{1}{p_{1}}$. There exists a constant $C$ depending continuously on $p,p_1$, such that
$$
\left\|\left(2^{j}\left\|[\Delta_{j},v \pa_x] f\right\|_{L^{p}}\right)_{j}\right\|_{\ell^{r}} \leq C\left(\|\pa_x v\|_{L^{\infty}}\|f\|_{B_{p, r}^{1}}+\|\pa_x f\|_{L^{p_{2}}}\|\pa_x v\|_{B_{p_1,r}^{0}}\right).
$$
\end{lemma}
Finally, let us end this section with the key estimates.
\begin{lemma}[see Lemma 3.26 in \cite{B}]\label{lem2.3}
Let $u_0\in C^{0,1}(\R)$.
Assume that $\pa_xu\in L^1([0,T]; L^\infty(\R))$ solves \eqref{CH}.
Then there exists a universal constant $C$ such that for all $t\in[0,T]$, we have
\begin{align*}
&\|u(t)\|_{C^{0,1}(\R)}\leq \|u_0\|_{C^{0,1}(\R)}\exp\left(C\int_0^t \|\pa_xu(\tau)\|_{L^{\infty}(\R)}\mathrm{d}\tau\right).
\end{align*}
\end{lemma}

\section{Proof of Theorem \ref{th1}}\label{sec3}
In this section, we shall present details for the proof of the main theorem.
\subsection{Construction of Initial Data}\label{sec3.1}

Define a smooth cut-off function $\chi$ with values in $[0,1]$ which satisfies
\bbal
\chi(\xi)=
\bca
1, \quad \mathrm{if} \ |\xi|\leq \frac14,\\
0, \quad \mathrm{if} \ |\xi|\geq \frac12.
\eca
\end{align*}
From now on, we set $\gamma:=\frac{17}{24}$ just for the sake of simplicity. Letting $$n\in 16\mathbb{N}=\left\{16,32,48,\cdots\right\}$$
 and
  $$\mathbb{N}(n)=\left\{k\in 8\mathbb{N}: \frac{n}4 \leq k\leq \frac{n}2\right\},$$
we can define the initial data $u_0$ by
\bbal
u_0(x)&=2^{-n}n^{-\frac{1}{2}}\log n\sum_{\ell\in \mathbb{N}(n)}\left(\cos(x+2^{\ell+1}\gamma)\cd(2^n\gamma+2^{\ell}\gamma)
+\cos(x+2^{\ell+1}\gamma)\cd(2^n\gamma-2^{\ell}\gamma)\right)
\check{\chi}(x+2^{\ell+1}\gamma).
\end{align*}
\begin{lemma}\label{le-e1} Let $p\in(2,\infty]$.
There exists a positive constant $C$ independent of $n$ such that
\bbal
&\|u_0\|_{L^\infty}\leq C 2^{-n}n^{-\frac12}\log n,\\
&\|\pa_xu_0\|_{L^\infty}\leq C n^{-\frac12}\log n,\\
&\|u_0\|_{B^1_{p,1}}\leq C n^{-\frac12(1-\frac2p)}\log n.
\end{align*}
\end{lemma}
{\bf Proof.} Notice that the simple fact $\cos(a)=\fr12\left(\cos(a)+\cos(-a)\right)=\fr12\big(e^{ia}+e^{-ia}\big)$, then by easy computations, we give that
\bbal
u_0&=2^{-n-1}n^{-\frac{1}{2}}\log n\sum_{\ell\in \mathbb{N}(n)}\left({\Phi^{++}_\ell}+{\Phi^{+-}_\ell}
+{\Phi^{-+}_\ell}+{\Phi^{--}_\ell}\right),
\end{align*}
where
\bbal
&\Phi^{++}_\ell=e^{i(x+2^{\ell+1}\gamma)\cd(2^n\gamma+2^{\ell}\gamma)}\check{\chi}(x+2^{\ell+1}\gamma), \\ &\Phi^{+-}_\ell=e^{i(x+2^{\ell+1}\gamma)\cd(2^n\gamma-2^{\ell}\gamma)}\check{\chi}(x+2^{\ell+1}\gamma),\\
&\Phi^{-+}_\ell=e^{i(x+2^{\ell+1}\gamma)\cd(-2^n\gamma+2^{\ell}\gamma)}\check{\chi}(x+2^{\ell+1}\gamma), \\ &\Phi^{--}_\ell=e^{i(x+2^{\ell+1}\gamma)\cd(-2^n\gamma-2^{\ell}\gamma)}\check{\chi}(x+2^{\ell+1}\gamma).
\end{align*}
It is easy to see that its Fourier transform
\bbal
\widehat{u_0}(\xi)=2^{-n-1}n^{-\frac{1}{2}}\log n\sum_{\ell\in \mathbb{N}(n)}\left(\widehat{\Phi^{++}_\ell}+\widehat{\Phi^{+-}_\ell}
+\widehat{\Phi^{-+}_\ell}+\widehat{\Phi^{--}_\ell}\right),
\end{align*}
where
\bbal
&\widehat{\Phi^{++}_\ell}=e^{i2^{\ell+1}\xi\cd \gamma}\chi(\xi-2^n\gamma-2^\ell\gamma), \\
&\widehat{\Phi^{+-}_\ell}=e^{i2^{\ell+1}\xi\cd \gamma}\chi(\xi-2^n\gamma+2^\ell\gamma),
\\
&\widehat{\Phi^{-+}_\ell}=e^{i2^{\ell+1}\xi\cd \gamma}\chi(\xi+2^n\gamma-2^\ell\gamma), \\
&\widehat{\Phi^{--}_\ell}=e^{i2^{\ell+1}\xi\cd \gamma}\chi(\xi+2^n\gamma+2^\ell\gamma),
\end{align*}
which implies that
\bbal
\mathrm{supp} \ \widehat{u_0}\subset \left\{\xi\in\R: \ 2^{n}\gamma-2^{\ell}\gamma-\fr12\leq |\xi|\leq 2^{n}\gamma+2^{\ell}\gamma+\fr12\right\}.
\end{align*}
Since $\check{\chi}$ is a Schwartz function, we have
\bbal
|\check{\chi}(x)|\leq C(1+|x|)^{-M}, \qquad  M\gg 1,
\end{align*}
then
\bbal
\|u_0\|_{L^\infty}&\leq C2^{-n}n^{-\frac{1}{2}}\log n\left\|\sum_{\ell\in \mathbb{N}(n)}\check{\chi}(x+2^{\ell+1}\gamma)\right\|_{L^\infty}
\\&\leq C2^{-n}n^{-\frac{1}{2}}\log n\left\|\sum_{\ell\in \mathbb{N}(n)}\frac{1}{(1+|x+2^{\ell+1}\gamma|)^M}\right\|_{L^\infty}\\
&\leq C2^{-n}n^{-\frac{1}{2}}\log n,
\end{align*}
and
\bbal
\|u_0\|_{L^1}&\leq C2^{-n}n^{-\frac{1}{2}}\log n\sum_{\ell\in \mathbb{N}(n)}\int_{\R}|\check{\chi}(x+2^{\ell+1}\gamma)|\dd x
\\&\leq C2^{-n}n^{-\frac{1}{2}}\log n\sum_{\ell\in \mathbb{N}(n)}\int_{\R}\frac{1}{(1+|x+2^{\ell+1}\gamma|)^M}\dd x\\
&\leq C2^{-n}n^{\frac12}\log n.
\end{align*}
Due to the fact $\Delta_ju_0=u_0$ if $j=n-1$ and $\Delta_ju_0=0$ if $j\neq n-1$, using the classical interpolation inequality $L^p=[L^1,L^\infty]_{(1/p,(p-1)/p)}$, we deduce that
\bbal
\|u_0\|_{B^1_{p,1}}&\leq C2^n\|u_0\|_{L^p}\leq C n^{-\frac12(1-\frac2p)}\log n.
\end{align*}
This completes the proof of Lemma \ref{le-e1}.

\begin{lemma}\label{le-e2} Let $p\in(2,\infty]$.
There exists a positive constant $c$ independent of $n$ such that
\bbal
\left\|(\pa_xu_0)^2\right\|_{B^0_{p,1}\left(\mathbb{N}(n)\right)}\geq c\log^2n, \qquad n\gg1.
\end{align*}
\end{lemma}
{\bf Proof.} Since $\mathrm{supp }\ \chi(\cdot-a)*\chi(\cdot-b)\subset B(a+b,1)$, we see that
\bbal
\mathrm{supp} \ \widehat{\Phi^{\nu\lambda}_\ell}*\widehat{\Phi^{\nu\mu}_m}\subset B(\nu2^{n+1}\gamma,2^{2+\frac n2}), \quad \lambda,\mu,\nu\in\{+,-\}.
\end{align*}
Moreover, noticing that
 $$\mathrm{supp} \ \widehat{\Phi^{++}_\ell}*\widehat{\Phi^{--}_m}\subset B((2^\ell-2^m)\gamma,1),$$
  we see that $$\mathrm{supp} \ \widehat{\Phi^{++}_\ell}*\widehat{\Phi^{--}_\ell}\subset B(0,1).$$
It follows that for any $j\in \mathbb{N}(n)$,
\bal\label{l}
\De_j[(\pa_x u_0)^2]=2^{-2n-1}n^{-1}\log^2n
\cdot\De_jU.
\end{align}
where $$U:=\sum_{\ell,m\in \mathbb{N}(n)}\Big(\pa_x\Phi^{++}_\ell\pa_x\Phi^{-+}_m+\pa_x\Phi^{++}_\ell\pa_x\Phi^{--}_m
 +\pa_x\Phi^{+-}_\ell\pa_x\Phi^{-+}_m
+\pa_x\Phi^{+-}_\ell\pa_x\Phi^{--}_m\Big).$$
We rewrite $U$ as
\bbal
U&=\sum_{\ell\in \mathbb{N}(n)}\Big(\pa_x\Phi^{++}_\ell\pa_x\Phi^{-+}_\ell+\pa_x\Phi^{++}_\ell\pa_x\Phi^{--}_\ell
 +\pa_x\Phi^{+-}_\ell\pa_x\Phi^{-+}_\ell
+\pa_x\Phi^{+-}_\ell\pa_x\Phi^{--}_\ell\Big)\\
&\quad+\sum_{m\in \mathbb{N}(n)}\sum_{\ell<m\in \mathbb{N}(n)}\Big(\pa_x\Phi^{++}_\ell\pa_x\Phi^{-+}_m+\pa_x\Phi^{++}_\ell\pa_x\Phi^{--}_m
 +\pa_x\Phi^{+-}_\ell\pa_x\Phi^{-+}_m
+\pa_x\Phi^{+-}_\ell\pa_x\Phi^{--}_m\Big)\\
&\quad+\sum_{m\in \mathbb{N}(n)}\sum_{\ell>m\in \mathbb{N}(n)}\Big(\pa_x\Phi^{++}_\ell\pa_x\Phi^{-+}_m+\pa_x\Phi^{++}_\ell\pa_x\Phi^{--}_m
 +\pa_x\Phi^{+-}_\ell\pa_x\Phi^{-+}_m
+\pa_x\Phi^{+-}_\ell\pa_x\Phi^{--}_m\Big)\\
&:=U_1+U_2+U_3.
\end{align*}
To deal with term $U_1$, using the facts that
\bbal
&\mathrm{supp} \ \widehat{\Phi^{++}_j}*\widehat{\Phi^{--}_j}\subset B(0,1),\\
&\mathrm{supp} \ \widehat{\Phi^{+-}_j}*\widehat{\Phi^{-+}_j}\subset B(0,1),
\end{align*}
thus, for any $j\in \mathbb{N}(n)$,
\bal\label{l1}
\De_j(\pa_x\Phi^{++}_j\pa_x\Phi^{--}_j+\pa_x\Phi^{+-}_j\pa_x\Phi^{-+}_j)=0.
\end{align}
Note that
\bbal
&\mathrm{supp} \ \widehat{\Phi^{++}_j}*\widehat{\Phi^{-+}_j}\subset B(2^{j+1}\gamma,1)\subset 2^j\mathcal{C}\left(0,4/3,3/2\right),
\\&\mathrm{supp} \ \widehat{\Phi^{+-}_j}*\widehat{\Phi^{--}_j}\subset B(2^{j+1}\gamma,1)\subset 2^j\mathcal{C}\left(0,4/3,3/2\right),
\end{align*}
and
$
\varphi(2^{-j}\xi)\equiv 1$ for $\xi\in 2^j\mathcal{C}\left(0,4/3,3/2\right),
$
we have
\bal\label{l2}
\De_j(\pa_x\Phi^{++}_j\pa_x\Phi^{-+}_j+\pa_x\Phi^{+-}_j\pa_x\Phi^{--}_j)=\pa_x\Phi^{++}_j\pa_x\Phi^{-+}_j+\pa_x\Phi^{+-}_j\pa_x\Phi^{--}_j.
\end{align}
\eqref{l1} and \eqref{l2} gives that
\bal\label{l3}
\De_jU_1=\pa_x\Phi^{++}_j\pa_x\Phi^{-+}_j+\pa_x\Phi^{+-}_j\pa_x\Phi^{--}_j.
\end{align}
To deal with term $U_2$,  we have if $j<m\in \mathbb{N}(n)$
\bbal
&\mathrm{supp} \ \widehat{\Phi^{++}_\ell}*\widehat{\Phi^{-+}_m}\subset B((2^\ell+2^m)\gamma,1)\subset 2^{j}\mathcal{C}(0,3,2^N),
\\&\mathrm{supp} \ \widehat{\Phi^{++}_\ell}*\widehat{\Phi^{--}_m}\subset B((2^\ell-2^m)\gamma,1)\subset 2^{j}\mathcal{C}(0,3,2^N),
\\&\mathrm{supp} \ \widehat{\Phi^{+-}_\ell}*\widehat{\Phi^{-+}_m}\subset B((-2^\ell+2^m)\gamma,1)\subset 2^{j}\mathcal{C}(0,3,2^N),
\\&\mathrm{supp} \ \widehat{\Phi^{+-}_\ell}*\widehat{\Phi^{--}_m}\subset B((-2^\ell-2^m)\gamma,1)\subset 2^{j}\mathcal{C}(0,3,2^N),
\end{align*}
and if $j>m\in \mathbb{N}(n)$
\bbal
&\mathrm{supp} \ \widehat{\Phi^{++}_\ell}*\widehat{\Phi^{-+}_m}\subset B((2^\ell+2^m)\gamma,1)\subset 2^{j}B(0,1/2),
\\&\mathrm{supp} \ \widehat{\Phi^{++}_\ell}*\widehat{\Phi^{--}_m}\subset B((2^\ell-2^m)\gamma,1)\subset 2^{j}B(0,1/2),
\\&\mathrm{supp} \ \widehat{\Phi^{+-}_\ell}*\widehat{\Phi^{-+}_m}\subset B((-2^\ell+2^m)\gamma,1)\subset 2^{j}B(0,1/2),
\\&\mathrm{supp} \ \widehat{\Phi^{+-}_\ell}*\widehat{\Phi^{--}_m}\subset B((-2^\ell-2^m)\gamma,1)\subset 2^{j}B(0,1/2).
\end{align*}
Thus,
\bal\label{l4}
\De_jU_2&=\sum_{j>\ell\in \mathbb{N}(n)}\Big(\pa_x\Phi^{++}_\ell\pa_x\Phi^{-+}_j+\pa_x\Phi^{++}_\ell\pa_x\Phi^{--}_j
 +\pa_x\Phi^{+-}_\ell\pa_x\Phi^{-+}_j
+\pa_x\Phi^{+-}_\ell\pa_x\Phi^{--}_j\Big).
\end{align}
Similarly,
\bal\label{l5}
\De_jU_3&=\sum_{j>\ell\in \mathbb{N}(n)}\Big(\pa_x\Phi^{++}_j\pa_x\Phi^{-+}_\ell+\pa_x\Phi^{++}_j\pa_x\Phi^{--}_\ell
 +\pa_x\Phi^{+-}_j\pa_x\Phi^{-+}_\ell
+\pa_x\Phi^{+-}_j\pa_x\Phi^{--}_\ell\Big).
\end{align}
Therefore, plugging \eqref{l3}-\eqref{l5} into \eqref{l}, we have
$$
\De_j[(\pa_x u_0)^2]:=I_1+I_2+I_3+I_4+I_5
$$
where
$$
I_1=2^{-2n-1}n^{-1}\log^2n\cdot(\pa_x\Phi^{++}_j\pa_x\Phi^{-+}_j+\pa_x\Phi^{+-}_j\pa_x\Phi^{--}_j)
$$
and
$$I_2=2^{-2n-1}n^{-1}\log^2n\cdot\sum_{j>\ell\in \mathbb{N}(n)}(\pa_x\Phi^{++}_j\pa_x\Phi^{-+}_\ell+\pa_x\Phi^{++}_\ell\pa_x\Phi^{-+}_j)$$
$$I_3=2^{-2n-1}n^{-1}\log^2n\cdot\sum_{j>\ell\in \mathbb{N}(n)}(\pa_x\Phi^{++}_j\pa_x\Phi^{--}_\ell+\pa_x\Phi^{++}_\ell\pa_x\Phi^{--}_j)$$
$$I_4=2^{-2n-1}n^{-1}\log^2n\cdot\sum_{j>\ell\in \mathbb{N}(n)}(\pa_x\Phi^{+-}_j\pa_x\Phi^{-+}_\ell+\pa_x\Phi^{+-}_\ell\pa_x\Phi^{-+}_j)$$
$$I_5=2^{-2n-1}n^{-1}\log^2n\cdot\sum_{j>\ell\in \mathbb{N}(n)}(\pa_x\Phi^{+-}_j\pa_x\Phi^{--}_\ell+\pa_x\Phi^{+-}_\ell\pa_x\Phi^{--}_j).$$
We will show that the term $I_1$ contributes the main part and the other terms are much less than the first one.
By direct computations give that
\bbal
I_1&=n^{-1}\log^2n\left(I_{1,1}+I_{1,2}+I_{1,3}\right).
\end{align*}
where
\bbal
&I_{1,1}=(1-2^{2j-2n})\gamma^2\cdot\cos[(x+2^{j+1}\gamma)\cd 2^{j+1}\gamma]\check{\chi}^2(x+2^{j+1}\gamma),\\
&I_{1,2}=2^{-2n}\cdot\cos[(x+2^{j+1}\gamma)\cd 2^{j+1}\gamma][\pa_x\check{\chi}(x+2^{j+1}\gamma)]^2,\\
&I_{1,3}=-2^{j-2n}\gamma\cdot\sin[(x+2^{j+1}\gamma)\cd 2^{j+1}\gamma][(\check{\chi}\pa_x\check{\chi})(x+2^{j+1}\gamma)].
\end{align*}
By Riemann's lemma, we have for large enough $n$,
\bbal
&\|I_{1,1}\|_{L^p}=(1-2^{2j-2n})\gamma^2\left\|\cos\left((2^{j+1}\gamma)x\right)\check{\chi}^2(x)\right\|_{L^p}\geq C,\\
&\|I_{1,2}\|_{L^p}\leq C2^{-2n},\\
&\|I_{1,3}\|_{L^p}\leq C2^{j-2n}\leq C2^{-n}.
\end{align*}
Thus, for large enough $n$, we have
\bal\label{y1}
\|I_1\|_{L^p}\geq cn^{-1}\log^2n.
\end{align}
Since $\mathrm{supp} \ \Phi^{\mu,\nu}_\ell\subset 2^n\mathcal{C},\ell\in \mathbb{N}(n), \mu,\nu\in\{+,-\}$, by Bernstein's inequality, we have
\bbal
\|I_2\|_{L^\infty}&\leq Cn^{-1}\log^2n2^{-2n}\left\|\sum_{j>\ell\in \mathbb{N}(n)}(\pa_x\Phi^{++}_j\pa_x\Phi^{-+}_\ell+\pa_x\Phi^{++}_\ell\pa_x\Phi^{-+}_j)\right\|_{L^\infty}
\\&\leq Cn^{-1}\log^2n\sum_{j>\ell\in \mathbb{N}(n)}\left\|(1+|x+2^{j+1}\gamma|^2)^{-M}(1+|x+2^{\ell+1}\gamma|^2)^{-M}\right\|_{L^\infty}\\&\leq Cn^{-1}\log^2n\sum_{j>\ell\in \mathbb{N}(n)}\left\|(1+|x|^2)^{-M}(1+|x-(2^{j+1}-2^{\ell+1})\gamma|^2))^{-M}\right\|_{L^\infty}\\&\leq Cn^{-1}\log^2n\sum_{j>\ell\in \mathbb{N}(n)}\left(\gamma (2^{j}-2^{\ell})\right)^{-2M}
\\&\leq C\log^2n2^{-\frac{nM}{2}},
\end{align*}
where we have separated $\R$ into two different regions $\{x: |x|\leq\gamma (2^{j}-2^{\ell})\}$ and $\{x: |x| > \gamma (2^{j}-2^{\ell})\}$.

We easily see that
\bbal
\|I_2\|_{L^1}\leq C\log^2n.
\end{align*}
Using the classical interpolation inequality leads to
\bal\label{y2}
\|I_2\|_{L^p}\leq \|I_2\|^{\frac{1}{p}}_{L^1}\|I_2\|^{\frac{p-1}{p}}_{L^\infty}\leq C\log^2n2^{-\frac{1}{2}(1-\frac{1}{p})nM}.
\end{align}
Applying similar argument as in \eqref{y2}, we have
\bal\label{y3}
\|I_i\|_{L^p}\leq C\log^2n2^{-\frac{1}{2}(1-\frac{1}{p})nM} \quad \mathrm{for} \quad i=3,4,5.
\end{align}
Combining \eqref{y1}-\eqref{y3}, we deduce that for $n\gg1$,
\bbal
\left\|\De_j[(\pa_x u_0)^2]\right\|_{L^p}\geq cn^{-1}\log^2n.
\end{align*}
Therefore, by the definition of the Besov norm, we have
\bbal
\left\|(\pa_xu_0)^2\right\|_{B^0_{p,1}(\mathbb{N}(n))}\geq \sum_{j\in \mathbb{N}(n)}\left\|\De_j[(\pa_xu_0)^2]\right\|_{L^p}\geq c\log^2n.
\end{align*}
This completes the proof of lemma \ref{le-e2}.

\subsection{The Equation Along the Flow}\label{sec3.2}

Given a Lipschitz velocity field $u$, we may solve the following ODE to find the flow induced by $u$:
\begin{align}\label{ode}
\quad\begin{cases}
\frac{\dd}{\dd t}\phi(t,x)=u(t,\phi(t,x)),\\
\phi(0,x)=x,
\end{cases}
\end{align}
which is equivalent to the integral form
\bal\label{n}
\phi(t,x)=x+\int^t_0u(\tau,\phi(\tau,x))\dd \tau.
\end{align}
Considering the equation
\begin{align}\label{pde}
\quad\begin{cases}
\pa_tv+u\pa_xv=P,\\
v(0,x)=v_0(x),
\end{cases}
\end{align}
then, we get from \eqref{pde} that
\bbal
\pa_t(\De_jv)+u\pa_x\De_jv&=R_j+\Delta_jP,
\end{align*}
with $R_j=[u,\De_j]\pa_xv=u\De_j\pa_xv-\Delta_j(u\pa_xv)$.

Due to \eqref{ode}, then
\bbal
\frac{\dd}{\dd t}\left((\De_jv)\circ\phi\right)&=R_j\circ\phi+\Delta_jP\circ\phi,
\end{align*}
which means that
\bal\label{l6}
\De_jv\circ\phi=\De_jv_0+\int^t_0R_j\circ\phi\dd \tau+\int^t_0\Delta_jP\circ\phi\dd \tau.
\end{align}
\subsection{Norm Inflation}\label{sec3.3}

For $n\gg1$, using Lemma \ref{lem2.1}, we have for $t\in[0,1]$
\bbal
\|u\|_{C^{0,1}}\leq C\|u_0\|_{C^{0,1}}\leq Cn^{-\frac12}\log n,
\end{align*}
and also
\bbal
\|u\|_{W^{1,p}}\leq C\|u_0\|_{W^{1,p}}\leq Cn^{-\frac12(1-\frac2p)}\log n.
\end{align*}
To prove Theorem \ref{th1}, it suffices to show that there exists $t_0\in(0,\frac{1}{\log n}]$ such that
\bal\label{m}
\|u(t_0,\cdot)\|_{B^1_{p,1}}\geq \log\log n.
\end{align}
We prove \eqref{m} by contraction. If \eqref{m} were not true, then
\bal\label{m1}
\sup_{t\in(0,\frac{1}{\log n}]}\|u(t,\cdot)\|_{B^1_{p,1}}< \log\log n.
\end{align}
We divide the proof into two steps.

{\bf Step 1: Lower bounds for $(\De_ju)\circ \phi$}

Now we consider the equation along the Lagrangian flow-map associated to $u$.
Utilizing \eqref{l6} to \eqref{CH} yields
\bbal
(\De_ju)\circ \phi&=\De_ju_0+\int^t_0R^1_j\circ \phi\dd \tau +\int^t_0\De_jF\circ \phi\dd \tau
\\&\quad +\int^t_0\big(\De_jE\circ \phi-\De_jE_0\big)\dd \tau+t\De_jE_0,
\end{align*}
where
\bbal
&R^1_j=[u,\De_j]\pa_xu,
\\&F=-\pa_x(1-\pa^2_x)^{-1}u^2, \\
&E=-\frac12\pa_x(1-\pa^2_x)^{-1}(\pa_xu)^2.
\end{align*}
Due to Lemma \ref{le-e2}, we deduce
\bal\label{g1}
\sum_{j\in \mathbb{N}(n)}2^j\|\De_jE_0\|_{L^p}
&\approx \sum_{j\in \mathbb{N}(n)}\|\De_j\pa_xE_0\|_{L^p}\nonumber\\
&\geq c\sum_{j\in \mathbb{N}(n)}\|\De_j(\pa_xu_0)^2\|_{L^p}\nonumber\\
&\geq c\log^2n.
\end{align}
Notice that \eqref{n}, then we have for $t\in(0,\frac{1}{\log n}]$,
\bbal
\frac12\leq |\pa_x\phi|\leq 2,
\end{align*}
thus for $p\in(2,\infty)$
\bbal
\|f(t,\phi(t,\cdot))\|_{L^p}\approx_p \|f(t,\cdot)\|_{L^p}\quad\text{and}\quad\|f(t,\phi(t,x))\|_{L^\infty}= \|f(t,\cdot)\|_{L^\infty}.
\end{align*}
Then, using the commutator estimate from Lemma \ref{lem2.2}, we have
\bal\label{g2}
\sum_{j\geq -1}2^j\|R^1_j\circ \phi\|_{L^p}&\leq C\sum_{j\geq -1}2^j\|R^1_j\|_{L^p}\nonumber\\
&\leq C\|\pa_xu\|_{L^\infty}\|u\|_{B^1_{p,1}}\nonumber\\
&\leq C n^{-\frac12}\log^3n.
\end{align}
Also, we have
\bbal
2^j\|\De_jF\circ \phi\|_{L^p}&\leq C2^j\|\De_jF\|_{L^p}
\\&\leq C\|u\|_{L^\infty}\|u\|_{L^p}
\\&\leq C(\|u_0\|_{L^\infty}+\|u_0\|^2_{W^{1,\infty}})\|u\|_{W^{1,p}}
\\&\leq Cn^{-1}n^{-\frac12(1-\frac2p)}\log^3n,
\end{align*}
which implies
\bal\label{g3}
\sum_{j\in \mathbb{N}(n)}2^j\|\De_jF\circ \phi\|_{L^p}\leq C n^{-\frac12(1-\frac2p)}\log^3n.
\end{align}
Combining \eqref{g1}-\eqref{g3} and using Lemmas \ref{le-e1}-\ref{le-e2} yields
\bbal
\sum_{j\in \mathbb{N}(n)}2^j\|(\De_ju)\circ \phi\|_{L^p}
&\geq t\sum_{j\in \mathbb{N}(n)}2^j\|\De_jE_0\|_{L^p}-\sum_{j\in \mathbb{N}(n)}2^j\|\De_jE\circ \phi-\De_jE_0\|_{L^p}-
Cn^{-\frac12+\frac1p}\log^3n
-C\|u_0\|_{B^1_{p,1}}
\\&\geq ct\log^2n-\sum_{j\in \mathbb{N}(n)}2^j\|\De_jE\circ \phi-\De_jE_0\|_{L^p}-
Cn^{-\frac12+\frac1p}\log n.
\end{align*}
{\bf Step 2: Upper bounds for $\De_jE\circ \phi-\De_jE_0$}

By easy computations, we can see that
\bal\label{E}
\pa_tE+u\pa_xE=G,
\end{align}
where
\bbal
G=&\frac13u^3-\frac12u(1-\pa^2_x)^{-1}(\pa_xu)^2
-(1-\pa^2_x)^{-1}\left(\frac13u^3-
\frac12u(\pa_xu)^2-\pa_x\Big(\pa_xu(1-\pa^2_x)^{-1}
(u^2+\frac12(\pa_xu)^2)\Big)\right).
\end{align*}
Utilizing \eqref{l6} to \eqref{E} yields
\bbal
\De_jE\circ \phi-\De_jE_0=\int^t_0[u,\De_j]\pa_xE\circ \phi\dd \tau +\int^t_0\De_jG\circ \phi\dd \tau.
\end{align*}
Using the commutator estimate from Lemma \ref{lem2.2}, one has
\bbal
2^j\|[u,\De_j]\pa_xE\|_{L^p}\leq C(\|\pa_xu\|_{L^\infty}\|E\|_{B^1_{p,\infty}}
+\|\pa_xE\|_{L^\infty}\|u\|_{B^1_{p,\infty}})\leq C\|u\|^2_{C^{0,1}}\|u\|_{W^{1,p}}
\end{align*}
and
\bbal
2^j\|\De_jG\|_{L^p}\leq C \|u\|^2_{C^{0,1}}\|u\|_{W^{1,p}}.
\end{align*}
Then, we have
\bbal
2^j\|\De_jE\circ \phi-\De_jE_0\|_{L^p}\leq \|u\|^2_{C^{0,1}}\|u\|_{W^{1,p}}\leq Cn^{-1}n^{-\frac12(1-\frac2p)}\log^3n,
\end{align*}
which leads to
\bbal
\sum_{j\in \mathbb{N}(n)}2^j\|\De_jE\circ \phi-\De_jE_0\|_{L^p}
\leq Cn^{-\frac12+\frac1p}\log^3n.
\end{align*}
Combining Step 1 and Step 2, then for $t=\frac{1}{\log n}$, we obtain for $n\gg1$
\bbal
\|u(t)\|_{B^1_{p,1}}&\geq \|u(t)\|_{B^1_{p,1}(\mathbb{N}(n))}\\
&\geq C\sum_{j\in \mathbb{N}(n)}2^j\|(\De_ju)\circ \phi\|_{L^p}
\\&\geq c t\log^2n-Cn^{-\frac12+\frac1p}\log n\\
&\geq \log\log n,
\end{align*}
which contradicts the hypothesis \eqref{m1}.

In conclusion, we obtain the norm inflation and hence the ill-posedness of the CH equation. Thus, Theorem \ref{th1} is proved.
\vspace*{1em}
\section*{Acknowledgements}
J. Li is supported by the National Natural Science Foundation of China (11801090 and 12161004) and Jiangxi Provincial Natural Science Foundation (20212BAB211004). Y. Yu is supported by the National Natural Science Foundation of China (12101011) and Natural Science Foundation of Anhui Province (1908085QA05). Y. Guo is supported by the Guangdong Basic and Applied Basic Research Foundation (No. 2020A1515111092). W. Zhu is partially supported by the National Natural Science Foundation of China (11901092) and Natural Science Foundation of Guangdong Province (2017A030310634).

\vspace*{1em}
\section*{Data Availability} Data sharing is not applicable to this article as no new data were created or analyzed in this study.

\section*{Conflict of interest}
The authors declare that they have no conflict of interest.

\end{document}